\documentclass[fleqn,11pt]{article}

\usepackage{amsmath, amscd, amsfonts, amssymb, graphicx, color}
\usepackage{hyperref}
\usepackage{cite}
\usepackage{amsmath}
\usepackage{amsfonts}
\usepackage{epsfig}
\usepackage{amssymb}
\usepackage{amsthm}
\usepackage{epstopdf}
\usepackage{footnote}
\usepackage{newlfont}
\usepackage{color}
\newtheorem{theorem}{Theorem}[section]
\newtheorem{lemma}[theorem]{Lemma}

\newtheorem{corollary}[theorem]{Corollary}
\theoremstyle{definition}
\newtheorem{definition}[theorem]{Definition}
\newtheorem{example}[theorem]{Example}
\theoremstyle{example}

\newtheorem{remark}[theorem]{Remark}

\begin{document}
%%%%%%%%%%%%%%%%%%%%%%%%%%%%%%%%%%%%%%%%%%%%%%%%%%%%%%%%

 \title{Nearest matrix with prescribed eigenvalues and its applications}

 \author{E. Kokabifar\thanks{Department of Mathematics,
 Faculty of Science, Yazd University, Yazd, Iran
 (e.kokabifar@stu.yazd.ac.ir, loghmani@yazd.ac.ir).},\,
 G.B. Loghmani\footnotemark[1]\,\,
 and S.M. Karbassi\thanks{Department of Mathematics, Yazd Branch, Islamic Azad
 University, Yazd, Iran, (mehdikarbassi@gmail.com).}}
\date{}
\maketitle

\vspace{-6mm}

\begin{abstract}
Consider $n \times n$ matrix $A$ and a set $\Lambda$ consisting of $k \le n$ prescribed complex numbers. Lippert (2010) in a challenging article, studied geometrically the spectral norm distance from $A$ to the set $\Lambda$ and constructed a perturbation matrix $\Delta$ with minimum spectral norm such that $A+\Delta$ had $\Lambda$ in its spectrum. This paper presents an easy practical computational method for constructing  the optimal perturbation $\Delta$ by extending necessary definitions and lemmas of previous works. Also, some conceivable applications of this issue are provided.

\end{abstract}



%%%%%%%%%%%%%%%%%%%%%%%%%%%%%%%%%%%%%%%%%%%%%%%%%%%%%%%%%%%%%%%%%%%%%%%
\section{Introduction}\label{intr}
Let $A$ be an $n\times n$ complex matrix and  let $L$ be the set of complex $n\times n$ matrices that have $\lambda \in \mathbb{C}$ as a prescribed multiple eigenvalue. In 1999, Malyshev \cite{malyshev} obtained the following formula for the spectral norm distance from $A$ to $L:$
\begin{equation*}
\mathop {\min }\limits_{B \in L} {\left\| {A - B} \right\|_2} = \mathop {\max }\limits_{\gamma  \ge 0} {s_{2n - 1}}\left( {\left[ {\begin{array}{*{20}{c}}
{A - \lambda I}&{\gamma {I_n}}\\
0&{A - \lambda I}
\end{array}} \right]} \right),
\end{equation*}
where $\|\cdot\|_2$ denotes the spectral matrix norm and $\,s_1(\cdot)
\geq s_2(\cdot) \geq s_3(\cdot) \geq \cdots \,$ are the singular
values of the corresponding matrix in nonincreasing order. Also he constructed a perturbation, $\Delta$, to matrix $A$ such that $A+\Delta$ belonged to the $L$ and $\Delta$ was the optimal perturbation of the matrix $A$. Malyshev's work can be considered as a solution to the Wilkinson's problem, that is, the computation of the distance from a matrix $A \in \mathbb{C}^{n \times n}$ which has only simple eigenvalues, to the set of $n \times n$ matrices with multiple eigenvalues. Wilkinson introduced this distance in \cite{wilkinson} and some bounds were computed for it by Ruhe \cite{ruhe}, Wilkinson \cite{wil1,wil2,wil3,wil4} and Demmel \cite{demmel1}. However, in a non-generic case, if $A$ is a normal matrix then the Malyshev's formula is not directly applicable. Ikramov and Nazari \cite{ikramovasli} showed this point and they obtained an extension of Malyshev's formula for normal matrices. Furthermore, the Malyshev's formula was extended by them \cite{ikramov} for the case of a spectral norm distance from $A$ to matrices with a prescribed triple eigenvalue. In 2011, under some conditions, a perturbation $\Delta$ to matrix $A$ was constructed by Mengi \cite{mengi} such that $\Delta$  had minimum spectral norm and $A+\Delta$ was belonged to the set of matrices that had a prescribed eigenvalue of prespecified algebraic multiplicity. Moreover, Malyshev's work also was extended by Lippert \cite{lipert} and Gracia \cite{gracia}. They computed a spectral norm distance from $A$ to the matrices with two prescribed eigenvalues. Recently, Lippert \cite{lipertk} introduced a geometric motivation for the results obtained in \cite{gracia, lipert} and he computed the smallest perturbation in the spectral norm such that the perturbed matrix had some given eigenvalues. On the other hand, in \cite{lipertk} and Section 6 of \cite{lipert} it was shown that the optimal perturbations are not always computable for the case of fixing three or more distinct eigenvalues. 
Denote by $\mathcal{M}_k$ the set of $n\times n$ matrices that have $k \le n$ prescribed eigenvalues. This article concerns the spectral norm distance from $A$ to $\mathcal{M}_k$ and describes a clear computational technique for construction of $\Delta$ having minimum spectral norm and satisfying $A+\Delta\in\mathcal{M}_k$. Also, some possible applications of this topic is considered. First, some lower bounds for this distance are obtained. Then, two assumptions are provided such that the optimal perturbation is always computable when these assumptions hold. It is noticeable that If one or both conditions are not satisfied then still $A+\Delta\in\mathcal{M}_k$, but $\Delta$ has not necessary minimum spectral norm. In this case we can have lower and upper bounds for the spectral norm distance form $A$ to $A+\Delta\in\mathcal{M}_k$. Note that if, in a special case, $A$ is a normal matrix, i.e., $A^*A=AA^*$, then we can not use the method described in this paper for the computation of the perturbation, immediately. In this case, by following the analysis performed in \cite{onaremarkable,ikramovasli,nazarirajabi} one can derive a refinement of our results for the case of normal matrices. Therefore, throughout of this paper, it is assumed that $A$ is not a normal matrix. Suppose now that an $n \times n$ matrix $A$ and a set of complex numbers $\Lambda  = \{ \lambda _1 ,\lambda _2 , \hdots ,\lambda _k \}$ in which $k \le n$, are given. Now for 
\[\gamma  = [{\gamma _{1,1}},{\gamma _{2,1}},{ { }} \ldots ,{\gamma _{k - 1,1}},{\gamma _{1,2}},{\gamma _{2,2}},{ { }} \ldots { {,}}{\gamma _{k - 2,2}},{ { }} \ldots ,{\gamma _{1,k - 1}}] \in {\mathbb{C}^{\frac{{k(k - 1)}}{2}}},\]
define  the $nk\times nk$ upper triangular matrix $Q_A(\gamma)$ as
\begin{equation}\label{Q}
{Q_A}(\gamma ) = {\left[ {\begin{array}{*{20}{c}}
{{B_1}}&{{\gamma _{1,1}}{I_n}}&{{\gamma _{1,2}}{I_n}}& \ldots &{{\gamma _{1,k - 1}}{I_n}}\\
0&{{B_2}}&{{\gamma _{2,1}}{I_n}}& \ldots &{{\gamma _{2,k - 2}}{I_n}}\\
 \vdots & \ddots &{{B_3}}& \ddots & \vdots \\
{}&{}&{}& \ddots &{{\gamma _{k - 1,1}}{I_n}}\\
0& \ldots &{}&0&{{B_k}}
\end{array}} \right]_{nk \times nk}},
\end{equation}
where $\gamma_{i,1}, (i=1, \hdots, k-1)$  are pure real variables and $B_j=A-\lambda_jI_n,~(j=1, \hdots, k).$

Clearly, $Q_A(\gamma)$ can be assumed as a matrix function of variables $\gamma_{i,j}$ such that $i=1\hdots k-1,~ j=1\hdots k-i.$
Hereafter, for the sake of simplicity, the positive integer $nk-(k-1)$ is denoted by $\kappa$.
Assume that the spectral norm distance from $A$ to $\mathcal{M}_k$ is denoted by $\rho_2 (A,\mathcal{M}_k)$, i.e,
\[{\rho _2}(A,\mathcal{M}_k) ={\left\| \Delta \right\|_2}= \mathop {\min }\limits_{M \in \mathcal{M}_k} {\left\| {A - M} \right\|_2},\]
in what follows, some lower bounds for the optimal perturbation $\Delta$ such that $A+\Delta \in \mathcal{M}_k$ are obtained.
\section{Lower bounds for the optimal perturbation}\label{lowe}
Let us begin by considering ${s_\kappa }\left( {Q_A(\gamma )} \right)$ which is the $\kappa$th singular value of $Q_A(\gamma)$. First, note that ${s_\kappa }\left( {Q_A(\gamma )} \right)$ is a continuous function of variable $\gamma$. Also, if we define the unitary matrix $U$ of the form
\[U = \left[ {\begin{array}{*{20}{c}}
I_n&{}&{}&{}&0\\
{}&{ - I_n}&{}&{}&{}\\
{}&{}&I_n&{}&{}\\
{}&{}&{}& \ddots &{}\\
0&{}&{}&{}&{{{( - 1)}^{k - 1}}I_n}
\end{array}} \right]_{k \times k},\]
where $I_n$ is $n \times n$ identity matrix, then it is straightforward to see that
\[UQ_A(\gamma )U^* = {\left[ {\begin{array}{*{20}{c}}
{{B_1}}&{{-\gamma _{1,1}}{I_n}}&{{\gamma _{1,2}}{I_n}}& \ldots &{{\gamma _{1,k - 1}}{I_n}}\\
0&{{B_2}}&{{-\gamma _{2,1}}{I_n}}& \ldots &{{\gamma _{2,k - 2}}{I_n}}\\
 \vdots & \ddots &{{B_3}}& \ddots & \vdots \\
{}&{}&{}& \ddots &{{-\gamma _{k - 1,1}}{I_n}}\\
0& \ldots &{}&0&{{B_k}}
\end{array}} \right]_{nk \times nk}}.\]

Since under a unitary transformation the singular values of a square matrix are invariant, it follows that ${s_\kappa }\left( {Q_A(\gamma )} \right)$ is an even function with respect to all its real variables $\gamma_{i,1}, (i=1, \ldots, k-1)$. Therefore, without loss of generality, in remainder of the paper, we assume that $\gamma_{i,1}\geq 0,$ for every $i=1,\ldots, k-1$.
\begin{lemma}\label{lem1}
If $A \in {\mathbb{C}^{n \times n}}$ has $\lambda _1 ,\lambda _2 , \ldots ,\lambda _k$ as some of its eigenvalues, then for all $\gamma  \in {\mathbb{C}^{\frac{{k(k - 1)}}{2}}},$ it holds that ${s_\kappa }\left( {Q_A(\gamma )} \right)=0.$
\end{lemma}
\textbf{Proof.} Suppose that $\lambda _1 ,\lambda _2 , \ldots ,\lambda _k$ are some of the eigenvalues of $A$ corresponding to the associated eigenvectors ${e_1},{e_2}, \ldots ,{e_k}$, respectively. Evidently,
\begin{equation}\label{ab}
Ae_i  = \lambda _i e_i, \qquad {  \mbox{and}} \qquad B_j e_i  = \left( {\lambda _i  - \lambda _j } \right)e_i ,\qquad {i,j = 1, \cdots ,k}.
\end{equation}

It is shown that $Q_A(\gamma)$ has $k$ linearly independent eigenvectors corresponding to zero as one of its eigenvalues. This implies that ${s_\kappa }\left( {Q_A(\gamma )} \right)=0.$
Two cases are now considered.

Case 1. Let $\lambda _1, \ldots ,\lambda _k$ be $k$ distinct eigenvalues of $Q_A(\gamma)$. Consequently, $ {{e_1}, \ldots ,{e_k}}$ are $k$ linearly independent eigenvectors. Consider now the $nk\times 1$ vectors ${v^1},{v^2}, \ldots ,{v^k}$  such that ${v^1}=\left[e_1,0,\ldots,0\right]^T $ and introduce the remaining vectors by the following formula
\begin{equation}\label{vjha}
v_j^m  = \left\{ {\begin{array}{*{20}c}
   0 & {j > m}  \\
   {e_m } & {j = m}  \\
   {\frac{1}{{\lambda _j  - \lambda _m }}\sum\limits_{p = 1}^{m - j} {\left( {\gamma _{j,p} v_{p + j}^m } \right)e_m } } & {j = m - 1, \ldots ,1}  \\
\end{array}} \right.,
\end{equation}
where $v_j^m$ denotes the $j$th component of $v^m$. The elements of each vector $v^m$ should be computed recursively, starting from the last component. Clearly, ${v^1},{v^2}, \ldots ,{v^k}$ are $k$ linearly independent vectors. Using (\ref{ab}) and (\ref{vjha}), straightforward calculation yields
\begin{eqnarray*}
\left( {Q_A \left( \gamma  \right)v^m } \right)_i  &=& B_i v_i^m  + \sum\limits_{p = 1}^{m - i} {\left( {\gamma _{i,p} v_{p + i}^m } \right)}  \\&=& \frac{1}{{\lambda _i  - \lambda _m }}\sum\limits_{p = 1}^{m - i} {\left( {\gamma _{i,p} v_{p + i}^m } \right)B_i e_m }  + \sum\limits_{p = 1}^{m - i} {\left( {\gamma _{i,p} v_{p + i}^m } \right)}  \\&=&  - \sum\limits_{p = 1}^{m - i} {\left( {\gamma _{i,p} v_{p + i}^m } \right)}  + \sum\limits_{p = 1}^{m - i} {\left( {\gamma _{i,p} v_{p + i}^m } \right)}  \\&=& 0.
\end{eqnarray*}
Thus, the vectors ${{v^1},{v^2}, \ldots ,{v^k}}$ satisfy $Q_A(\gamma)v^m=0, (m=1,\ldots,k)$. Consequently, the dimension of the null space of $Q_A(\gamma)$ is at least $\kappa$ which means that ${s_\kappa }\left( {Q_A(\gamma )} \right)=0.$

%Case 2. Suppose that $\left\{ {{e_1},{e_2}, \hdots ,{e_k}} \right\}$ is a set of dependent vectors.
%So, we can assume that $\lambda_i$ is a defective eigenvalue of algebraic multiplicity $l_i$ that has $\left\{ {{e_1},{e_2}, \hdots ,{e_{l_i}}} \right\}$ as a dependent set of its eigenvectors. Now, this set of vectors can be replaced with $\left\{ {{{\bar e}_1},{{\bar e}_2}, \hdots ,{{\bar e}_{l_i}}} \right\}$ of linearly independent generalized eigenvectors corresponding to $\lambda_i$ such that $\left\{ {{\bar e_1},{\bar e_2}, \ldots ,{\bar e_k}} \right\}$ be a set of independent vectors. The proof is completed by what has already been proven in Case 1. \qquad $\square$

Case 2. Assume that some of the complex numbers $\lambda _1 ,\lambda _2 , \ldots ,\lambda _k$ are equal. Without loss
of generality, we can assume that $\lambda_1$ is an eigenvalue of algebraic multiplicity $l$, i.e.,  $\lambda _1 =\lambda _2 = \ldots =\lambda _l$. In this case, a new method is provided for constructing the vectors associated with $\lambda_1$, i.e., $v^1,v^2,\ldots,v^l$. Obviously, it is enough to construct a set of linearly independent vectors $\{v^1,v^2,\ldots,v^l,v^{l+1},\ldots,v^k\}$, for which $Q_A(\gamma)v^m=0, (m=1,\ldots,l)$. To do this, let $e_1$ be the right eigenvector for $\lambda_1$ and let $e_1, \bar e_1, \bar e_2, \ldots, \bar e_{l-1}$ form a chain of generalized eigenvectors of length $l$ associated with $\lambda_1$. So, we have
\[
\left( {A - \lambda _1 I} \right)e_1  = 0,\hspace{.2cm}\left( {A - \lambda _1 I} \right){\bar e_1}  = e_1,\hspace{.6cm}{  \mbox{and}}\hspace{.6cm}\left( {A - \lambda _1 I} \right){\bar e_j}  = \bar e_{j-1},\hspace{.2cm}(j=2,\ldots,l).
\]

Define the operator $\mathcal{V}$ as follows:
\begin{equation}\label{vv}
\mathcal{V}[{e_1}] = {{\bar e}_1}\qquad {  \mbox{and}} \qquad \mathcal{V}[{{\bar e}_i}] = {{\bar e}_{i + 1}},\qquad i=1,\hdots,l-1.
\end{equation}
Now, let ${v^1}=\left[e_1,0,\ldots,0\right]^T$ and define the vectors $ {,{v^2}, \ldots ,{v^l}}$ by
\begin{equation*}
v^m_j= \left\{ {\begin{array}{*{20}c}
   0 & {j > m}  \\
   {e_1 } & {j = m}  \\
   { - \sum\limits_{p = 1}^{m - j} {\left( {\gamma _{j,p} {\cal V}\left[ {v_{p + j}^m} \right]} \right)} } & {j = m - 1, \ldots ,1}  \\
\end{array}} \right.,
\end{equation*}
where $v_j^m$ denotes the $j$th component of $v^m$ and analogous to Case 1 elements of each vector $v^m$, can now be computed by the above recursive relation, beginning from the last element. The remaining  vectors associated with simple eigenvalues are constructed by the method presented in Case 1. Linearity independence of vectors $v^i, (i=1,\ldots,l)$ is concluded by their definition. Thus, ${v^1}, \ldots ,{v^k}$ are linearly independent vectors. Also, calculations similar to what was performed in Case 1, conclude that $ {{v^1}, \ldots ,{v^l}}$ are the $l$ eigenvectors associated to zero as an eigenvalue of $Q_A(\gamma)$, i.e., $Q_A(\gamma)v^i=0,$ for every $i=1,\ldots,l$. $\qquad\square$

Following corollary can be concluded by similar calculations of the above lemma.
\begin{corollary}
Let $\lambda _1 ,\lambda _2 , \ldots ,\lambda _k$ be some 
eigenvalues of $A \in {\mathbb{C}^{n \times n}}$. Then for all $\gamma  \in {\mathbb{C}^{\frac{{k(k - 1)}}{2}}}$  we have
\begin{equation}\label{coro}
\begin{array}{*{20}c}
  \vspace{.3cm} {s_n \left( {A - \lambda _i I} \right) = 0,\qquad i = 1, \ldots ,k,}  \\
   \vspace{.3cm}  {s_{2n - 1} \left( {\left[ {\begin{array}{*{20}c}
   {A - \lambda _i I} & {\gamma _{1,1} I}  \\
   0 & {A - \lambda _j I}  \\
\end{array}} \right]} \right) = 0,\qquad i,j = 1, \ldots ,k,}  \\
 {s_{3n - 2} \left( {\left[ {\begin{array}{*{20}c}
   {A - \lambda _i I} & {\gamma _{1,1} I} & {\gamma _{1,2} I}  \\
   0 & {A - \lambda _j I} & {\gamma _{2,1} I}  \\
   0 & 0 & {A - \lambda _l I}  \\
\end{array}} \right]} \right) = 0,\qquad i,j,l = 1, \ldots ,k,}  \\
    \vdots   \\
   {s_\kappa  \left( {Q_A (\gamma )} \right) = 0.}  \\
\end{array}
\end{equation}
\end{corollary}
\begin{lemma}
Suppose that $A \in {\mathbb{C}^{n \times n}}$ has $\lambda _1 ,\lambda _2 , \ldots ,\lambda _k$ as some of its eigenvalues. If $\Delta$ is the minimum norm perturbation such that $A+\Delta \in \mathcal{M}_k$, then for every $\gamma$,
\[
{\left\| \Delta  \right\|_2  \ge s_\kappa  \left( {Q_A (\gamma )} \right)}.
\]
\end{lemma}
\textbf{Proof.} By Lemma \ref{lem1}, we know that ${s_\kappa }\left( {{Q_{A + \Delta }}(\gamma )} \right) = 0$. Applying the Weyl inequalities for singular values (for example, see Corollary 5.1 of \cite{demmel}) to the relation ${Q_{A + \Delta }}(\gamma ) = Q_A(\gamma)+I_{k}\otimes \Delta,$ yields
\[
s_\kappa  \left( {Q_A (\gamma )} \right)=\left| {s_\kappa  \left( {Q_A (\gamma )} \right) - s_\kappa  \left( {Q_{A + \Delta } (\gamma )} \right)} \right| \le \left\| {I_k  \otimes \Delta } \right\|_2  = \left\| \Delta  \right\|_2. \qquad \square 
\]
Similar calculations as performed above for the singular value $s_\kappa  \left( {Q_A (\gamma )}\right) $, can be considered for  remainder of the singular values appearing in the left hand side of the equations (\ref{coro}). Thus, the following relations are deduced:
\[
\begin{array}{*{20}c}
    \vspace{.3cm} {\left\| \Delta  \right\|_2  \ge s_n \left( {A - \lambda _i I} \right),\qquad i = 1, \ldots ,k,}  \\
  \vspace{.3cm}   {\left\| \Delta  \right\|_2  \ge s_{2n - 1} \left( {\left[ {\begin{array}{*{20}c}
   {A - \lambda _i I} & {\gamma _{1,1} I}  \\
   0 & {A - \lambda _j I}  \\
\end{array}} \right]} \right),\qquad i,j = 1, \ldots ,k,}  \\
   {\left\| \Delta  \right\|_2  \ge s_{3n - 2} \left( {\left[ {\begin{array}{*{20}c}
   {A - \lambda _i I} & {\gamma _{1,1} I} & {\gamma _{1,2} I}  \\
   0 & {A - \lambda _j I} & {\gamma _{2,1} I}  \\
   0 & 0 & {A - \lambda _l I}  \\
\end{array}} \right]} \right) ,\qquad i,j,l = 1, \ldots ,k,}  \\
    \vdots   \\
   {\left\| \Delta  \right\|_2  \ge s_\kappa  \left( {Q_A (\gamma )} \right).}  \\
\end{array}
\]
Next corollary gives the main result of this section.
\begin{corollary}\label{alfaha}
Assume that $A \in \mathbb{C}^{n \times n}$ and a set $\Lambda$ consisting of $k \le n$ complex numbers are given. Then for all $\gamma  \in {\mathbb{C}^{\frac{{k(k - 1)}}{2}}}$, the optimal perturbation $\Delta$, satisfies
\begin{equation}\label{alfa}
\left\| \Delta  \right\|_2  \ge \max \left\{ {\alpha _1 ,\alpha _2 , \ldots ,\alpha _k } \right\},
\end{equation}
where
\[
\begin{array}{*{20}c}
  \vspace{.3cm} {\alpha _1  = \max \left\{ {s_n \left( {A - \lambda _i I} \right),\qquad i = 1, \ldots ,k} \right\},}  \\
    \vspace{.3cm} {\alpha _2  = \max \left\{ {s_{2n - 1} \left( {\left[ {\begin{array}{*{20}c}
   {A - \lambda _i I} & {\gamma _{1,1} I}  \\
   0 & {A - \lambda _j I}  \\
\end{array}} \right]} \right),\qquad i,j = 1, \ldots ,k} \right\},}  \\
   {\alpha _3  = \max \left\{ {s_{3n - 2} \left( {\left[ {\begin{array}{*{20}c}
   {A - \lambda _i I} & {\gamma _{1,1} I} & {\gamma _{1,2} I}  \\
   0 & {A - \lambda _j I} & {\gamma _{2,1} I}  \\
   0 & 0 & {A - \lambda _l I}  \\
\end{array}} \right]} \right) = 0,\hspace{.635cm} i,j,l = 1, \ldots ,k} \right\},}  \\
    \vdots   \\
   {\alpha _k  = s_\kappa  \left( {Q_A (\gamma )} \right).}  \\
\end{array}
\]
\end{corollary}
From this point, construction of an optimal perturbation is considered. This is completely dependent on dominance of $\alpha_i, (i=1,\ldots,k)$. First, let maximum of right hand side of (\ref{alfa}) occur in $\alpha_k={s_\kappa }\left( {Q_A(\gamma )} \right)$.
\section{Properties of ${s_\kappa }\left( {Q_A(\gamma )} \right)$ and its corresponding singular vectors }\label{prop}
In this section, we obtain further properties of ${s_\kappa }\left( {Q_A(\gamma )} \right)$ and its associated singular vectors. In the next section, this properties are applied to construct the optimal perturbation $\Delta$. For our discussion, it is necessary to reform some definitions and lemmas of \cite{malyshev, gracia, lipert,mengi}.
\begin{definition}\label{UV}
Suppose that vectors
\begin{equation*}
u(\gamma)=\left[ {\begin{array}{*{20}{c}}
{{u_1}(\gamma )}\\
 \vdots \\
{{u_k}(\gamma )}
\end{array}} \right], v(\gamma)=\left[ {\begin{array}{*{20}{c}}
{{v_1}(\gamma )}\\
 \vdots \\
{{v_k}(\gamma )}
\end{array}} \right] \in {\mathbb{C}^{nk}}~({u_j}(\gamma ),{v_j}(\gamma ) \in {\mathbb{C}^n},j = 1, \ldots ,k),
\end{equation*}
 is a pair of left and right singular vectors of ${s_\kappa }\left( {Q_A(\gamma )} \right)$, respectively.
Define
\begin{equation*}
U(\gamma) = [{u_1}(\gamma ), \ldots ,{u_k}(\gamma )]_{n \times k},\qquad {  \mbox{and}}\qquad V(\gamma) = [{v_1}(\gamma ), \ldots ,{v_k}(\gamma )]_{n \times k}.
\end{equation*}
\end{definition}

Considering definition of the vectors $u(\gamma)$ and $v(\gamma)$, we have the following relations
\begin{equation}\label{a}
Q_A(\gamma )v(\gamma ) = {s_\kappa }\left( {Q_A(\gamma )} \right)u(\gamma ),
\end{equation}\vspace{-.6cm}
\begin{equation}\label{b}
{Q_A}(\gamma )^*u(\gamma ) = {s_\kappa }\left( {Q_A(\gamma )} \right)v(\gamma ),
\end{equation}
also without loss of generality, assume that $u(\gamma)$ and $v(\gamma)$ are unit vectors.

The following lemma, which can be verified by considering Lemma 3.5 of \cite{lipert} (See also Lemma 3 of \cite{gracia}),  concludes that there exists a finite point $\gamma  \in {\mathbb{C}^{\frac{{k(k - 1)}}{2}}}$ where the function ${s_\kappa }\left( {Q_A(\gamma )} \right)$ attains its maximum value.
\begin{lemma}
${s_\kappa }\left( {Q_A(\gamma )} \right) \to 0$ as $\left| \gamma\right|  \to \infty$.
\end{lemma}
%\begin{lemma}
%If $\gamma_{i,j} \in \gamma$ for some $i,j $. Then ${s_\kappa }\left( {Q_A(\gamma )} \right) \to 0$ as $\left| \gamma_{i,j}\right|  \to \infty$.
%\end{lemma}
%\textbf{Proof.} It follows from Lemma 3.5 of \cite{lipert}. (See also Lemma 3 of \cite{gracia}).$\qquad \square$
\begin{definition}\label{sstar}
Let $\gamma_*$ be a point where the singular value ${s_\kappa }\left( {Q_A(\gamma )} \right)$ attains its maximum value such that $\mathop \Pi \limits_{i = 1}^{k - 1} {\gamma _*}_{i,1} > 0$ at this point. We set $\alpha_k^* ={s_\kappa }\left( {Q_A(\gamma_* )} \right).$
\end{definition}
It is easy to verify that if $\alpha_k^* =0$, then $\lambda_1, \lambda_2,\ldots,\lambda_k$ are some eigenvalues of $Q_A(\gamma)$. Therefore, in what follows we assume that $\alpha_k^*>0$. Moreover, suppose that $\alpha_k^*$ is a simple (not repeated) singular value of $Q_A(\gamma_*)$.

Investigating other properties of ${s_\kappa }\left( {Q_A(\gamma )} \right)$ and its associated  singular vectors requires results deduced in Theorem 2.10 and Theorem 2.11 of \cite{mengi}.
\begin{theorem}\label{meng}
Let $A(t)\colon \mathbb{R}\rightarrow \mathbb{C}^{n \times n}$ be an analytic matrix-valued function. There exists a decomposition $A(t)=U(t)S(t)V(t)^*$, where $U(t)\colon \mathbb{R}\rightarrow \mathbb{C}^{n \times n}, V(t)\colon \mathbb{R}\rightarrow \mathbb{C}^{n \times n}$ are unitary and analytic, $S(t)$ is diagonal and analytic for all $t$. Assume that $u_l(t), v_l(t)$ are the $l$th columns of $U(t)$ and $V(t)$, respectively, and $s_l(t)$ is the signed singular value at the $l$th diagonal entry of $S(t)$. Using the product rule and the fact that $A(t)v_l(t)=s_l(t)u_l(t)$, it is straightforward to deduce
\begin{equation}\label{mm1}
\frac{{ds_l (t)}}{{dt}} = {\mathop{  Real}} \left( {u_l (t)^*\frac{{dA(t)}}{{dt}}v_l (t)} \right).
\end{equation}
In particular, if $s_l(t_*)>0$ at a local extremum $t_*$, and $s_l(t_*)$ is a simple singular value, then 
\begin{equation}\label{mm2}
\frac{{ds_l (t_*)}}{{dt}} = {\mathop{  Real}} \left( {u_l (t_*)^*\frac{{dA(t_*)}}{{dt}}v_l (t_*)} \right)=0.
\end{equation}
\end{theorem}
Sun \cite{sun2} has shown that a simple singular value has an analytic expansion for the derivative formula similar to what is mentioned in (\ref{mm1}). In this case, we can differentiate a singular value along the tangent $\frac{{dA(t)}}{{dt}}$ from $A(t)$ and find a singular value decomposition that expresses the derivative as right hand side of (\ref{mm1}). The first part of the above theorem (Theorem 2.10 of \cite{mengi}) implies that for one parameter $A(t)$'s, Sun's formulas will still work, as long as $U$'s and $V$'s are chosen correctly. Moreover, the second part of Theorem \ref{meng} (Theorem 2.11 of \cite{mengi}) assures that singular vectors $u_l (t_*)$ and $v_l (t_*)$ for which equation (\ref{mm2}) is satisfied when $s_l(t_*)>0$, can be found. On the other hand, suppose we are given an $n$ parameter family $A(t_1,\ldots,t_n)$, and it is desired to find singular vectors $u$'s and $v$'s such that $\frac{{\partial s}}{{\partial t_i }} = {\mathop{  Real}} (u^* \frac{{\partial A}}{{\partial t_i }}v), (i=1,\ldots,n)$. In this case, Theorem \ref{meng} cannot be always applied. The $k$ eigenvalue case needs to deal with this failure analyticity. Fortunately, we can cope with this essential flaw by assuming that $\alpha_k^*$ is an isolated singular value (as this qualification is considered). Now, it should be noted that $\gamma_{i,1}, (i=1, \hdots, k-1)$  are pure real variables, while other  components of the vector $\gamma$ are complex variables. If we set ${\gamma _{i,j}} = {\gamma _{i,j,R}} + i{\gamma _{i,j,I}}$, for $j>1$, then it can be assumed that $\gamma$ has $(k-1)^2$ real components. Following lemma is deduced from applying Lemma \ref{meng} to $Q_A(\gamma).$

\begin{lemma}\label{lem3}
Let $\gamma_*$ and $\alpha_k^*$ be as defined in Definition \ref{sstar}, and let $\alpha_k^*>0$ be a simple singular value of $Q_A(\gamma_*)$. Then there exists a pair
\begin{equation*}
\left[ {\begin{array}{*{20}{c}}
{{u_1}({\gamma _*} )}\\
 \vdots \\
{{u_k}({\gamma _*} )}
\end{array}} \right],\left[ {\begin{array}{*{20}{c}}
{{v_1}({\gamma _*} )}\\
 \vdots \\
{{v_k}({\gamma _*} )}
\end{array}} \right] \in {\mathbb{C}^{nk}}({u_j}({\gamma _*} ),{v_j}({\gamma _*} ) \in {\mathbb{C}^n},j = 1, \ldots ,k),
\end{equation*}
of left and right singular vectors of $\alpha_k^*$, respectively, such that

$1.~u_i({\gamma _*} )^*{v_{i + j}}({\gamma _*} ) = 0,\qquad i = 1, \ldots k - 1,~j = 1, \ldots ,k-i,$ and

$2.~u_i({\gamma _*} )^*{u_j}({\gamma _*} ) = v_i({\gamma _*} )^*{v_j}({\gamma _*} ),\qquad i,j=1,\ldots,k$.
\end{lemma}
\textbf{Proof.} To prove first part of the lemma, two cases are considered.

\textit{Case }$1.$ At first assume that $j>1$.

By applying the Lemma 2.10 of \cite{mengi}, for both of the real and imaginary parts of ${{\gamma _*} _{i,j}} = {{\gamma _*} _{i,j,R}} + i{{\gamma _*} _{i,j,I}}$ we have Re$(u_i({\gamma _*} )^*{v_{i + j}}({\gamma _*} )) = 0,$ and Re$(i u_i({\gamma _*} )^*{v_{i + j}}({\gamma _*} )) = 0$, respectively. This means that $u_i({\gamma _*} )^*{v_{i + j}}({\gamma _*} ) = 0.$

\textit{Case }$2.$ Consider the case for which $j=1$.

Applying Lemma 2.10 of \cite{mengi} for real numbers ${\gamma _*}_{i,1}, (i=1,\hdots,k-1)$, concludes that Re$(u_i({\gamma _*} )^*{v_{i +1}}({\gamma _*} )) = 0$. So, first part of lemma is proved if we show that $u_i({\gamma _*} )^*{v_{i +1}}({\gamma _*} )$ is a real number. Let us introduce  the $1\times nk$ vectors $
\mathfrak u_{i,i}(\gamma_*)=[0, \hdots ,0,u_i({\gamma _*} )^*,0, \hdots ,0]$ and $\mathfrak v_{i,i}(\gamma_*)=[0, \hdots ,0,v_i({\gamma _*} )^*,0, \hdots ,0]$ where $u_i({\gamma _*} )^*$ and $v_i({\gamma _*} )^*$ are the $i$th components. By multiplying (\ref{a}) by $\mathfrak u_{i,i}(\gamma_*)$ and (\ref{b}) by $\mathfrak v_{i,i}(\gamma_*)$ both from the left, following results are obtained, respectively,
\begin{equation}\label{yek}
u_i({\gamma _*} )^*{B_i}{v_i}({\gamma _*} ) + \sum\limits_{p = 1}^{k - i} \left({{{\gamma _*} _{i,p}}u_i({\gamma _*} )^*{v_{i + p}}({\gamma _*} )}\right)  = \alpha_k^*u_i({\gamma _*} )^*{u_i}({\gamma _*} ),
\end{equation}\vspace{-.6cm}
\begin{equation}\label{dow}
v_i({\gamma _*} )^*B_i^*{u_i}({\gamma _*} ) + \sum\limits_{p = 1}^{i - 1} \left({{{\bar {\gamma _*} }_{p,i - p}}} v_i({\gamma _*} )^*{u_p}({\gamma _*} )\right) = \alpha_k^* v_i({\gamma _*} )^*{v_i}({\gamma _*} ),
\end{equation}
Conjugating (\ref{yek}), subtracting it from (\ref{dow}) and considering results obtained in Case 1 leads to the following relation
\begin{equation}\label{ghohar}
{{\gamma _*} _{i,1}}u_i({\gamma _*} )^*{v_{i + 1}}({\gamma _*} ) = \alpha_k^*\left( {u_i({\gamma _*} )^*{u_i}({\gamma _*} ) - v_i({\gamma _*} )^*{v_i}({\gamma _*} )} \right),
\end{equation}
clearly, the right hand side of (\ref{ghohar}) is a real number which implies that $u_i({\gamma _*} )^*{v_{i + 1}}({\gamma _*} )$ is also real. Thus firs part of the lemma is proved completely.

Now attempt to prove second part of the lemma. Multiplying (\ref{a}) by $\mathfrak u_{i,j}(\gamma_*)$ from the left where $u_i({\gamma _*} )^*$ is in the $j$th place, and multiplying (\ref{b}) by $\mathfrak v_{j,i}(\gamma_*)$ from the left where $v_j({\gamma_*})^*$ is in the $i$th place, leads to similar results as in (\ref{yek}) and (\ref{dow}), respectively. Performing similar calculations as first part of the proof, leads to
\[u_i({\gamma _*} )^*{B_j}{v_i}({\gamma _*} ) - u_i({\gamma _*} )^*{B_i}{v_j}({\gamma _*} ) = \alpha_k^* \left( {u_i({\gamma _*} )^*{u_i}({\gamma _*} ) - v_i({\gamma _*} )^*{v_i}({\gamma _*} )} \right),
\]
which can be written as
\[\left( {{\lambda _i} - {\lambda _j}} \right)u_i({\gamma _*} )^*{v_j}({\gamma _*} ) = \alpha_k^*\left( {u_i({\gamma _*} )^*{u_i}({\gamma _*} ) - v_i({\gamma _*} )^*{v_i}({\gamma _*} )} \right).\]
The left hand side of this equation is equal to zero. In fact, the assertion is obvious when $i=j.$ For $i\neq j$, assuming $i<j$ then considering first part of the proof. Thus proof is completed by keeping in mind that $\alpha_k^*>0$. \qquad $\square$
\begin{corollary}\label{natijeh}
Let suppositions of Lemma \ref{lem3} hold. Then the matrices $U({\gamma _*})$ and $V({\gamma _*})$ (recall Definition \ref{UV} but for $\gamma=\gamma_*$), satisfy ${U}({\gamma _*} )^*U({\gamma _*} ) = {V}({\gamma _*} )^*V({\gamma _*} )$.
\end{corollary}
\begin{lemma}\label{fullrank}
Suppose that assumptions of Lemma \ref{lem3} are satisfied. Then the two matrices $U({\gamma _*})$ and $V({\gamma _*})$ have full rank.
\end{lemma}
\textbf{Proof.} This result can be verified by considering Lemma 4.4 of \cite{lipert}. \qquad $\square$
%_______________________________________________________________
%Firstly suppose that maximum in right hand side of (\ref{max}) occurs in $\sigma _m, (1 \le m \le k).$ Let
%\[\Delta  =  - \left[ {{u_1}, \ldots ,{u_k}} \right]\left[ {\begin{array}{*{20}{c}}
%{{\sigma _1}}&{}&0\\
%{}& \ddots &{}\\
%0&{}&{{\sigma _k}}
%\end{array}} \right]{\left[ {{v_1}, \ldots ,{v_k}} \right]^\dag },\]
%where $\left[ {{v_1}, \ldots ,{v_k}} \right]^\dag $ is the Moore-Penrose pseudoinverse of $\left[ {{v_1}, \ldots ,{v_k}} \right]$. 
%
%Let $u_i ,v_i \in \mathbb{C}^n,~(i=1, \hdots, k)$ be a pair of left and right
%singular vectors of $B_i$, corresponding to
%$\sigma_i=s_n(B_i)$, respectively. It follows $B_iv_i=\sigma_i u_i,~(i=1, \hdots, k)$. Assume that singular vectors $v_1, \hdots ,v_k $ and $u_1, \hdots ,u_k $ satisfy
%\[\begin{array}{l}
%v_i^*{v_{i + j}} = 0,\qquad i = 1, \ldots ,k - 1,~j = 1, \ldots ,i - k,\\
%u_i^*{v_i} \ne 0,\qquad i = 1, \ldots ,k.
%\end{array}\]
% If we set 
% \[M = \sum\limits_{i = 1}^k {\left( {\frac{{{\sigma _i}}}{{\left\| {{v_i}} \right\|_2^2}}{u_i}v_i^*} \right)}, \]
%it is easy to show that $M{v_i} = {\sigma _i}{u_i},(i = 1, \ldots ,k)$. This implies that $\sigma_i, u_i$ and $v_i, (i=1, \hdots, k) $ are singular values and singular vectors of $M$. So, ${\left\| \Delta  \right\|_2}=\sigma_m.$ Moreover,
%\[\left( {{B_i} + \Delta } \right){v_i} = {\sigma _i}{u_i} - {\sigma _i}{u_i} = 0,\qquad (i = 1, \ldots ,k),\]
%or equivalently  
%\[\left( {A + \Delta } \right){v_i} = {\lambda _i}{v_i},\qquad(i = 1, \ldots ,k).\]
%________________________________________________________________
\section{Construction of the optimal perturbation }\label{const}
Let $\gamma_*$ and $\alpha_k^*$ be as defined in Definition \ref{sstar}, and let $\alpha_k^*>0$ be a simple singular value of $Q_A(\gamma_*)$. In this section, a perturbation $\Delta \in \mathbb{C}^{n \times n}$ is constructed such that satisfies ${\left\| \Delta  \right\|_2} =\alpha_k^*$, and perturbed matrix $A+\Delta$ has $k$ prescribed eigenvalues.
Suppose that the matrices $U(\gamma_* )$ and $ V(\gamma_* )$ are as in Corollary \ref{natijeh} and define
\begin{equation}\label{delta}
\Delta  =  - {\alpha_k^*}U({\gamma _*}){V }({\gamma _*})^\dag,
\end{equation}
where ${V }({\gamma _*})^\dag$ denotes the Moore-Penrose pseudoinverse of ${V }({\gamma _*})$. From Corollary \ref{natijeh} it can be deduced that the two matrices $U(\gamma _* )$ and $V(\gamma _* )$ have the same nonzero singular values. So, there exists a  unitary matrix $W \in \mathbb{C}^{n \times n}$ such that $U(\gamma_*)=W V(\gamma_*)$. Notice that Lemma \ref{fullrank} implies $ V{(\gamma )^\dag } V{(\gamma )}=I_k$. Therefore,
\begin{equation*}
{\left\| \Delta  \right\|_2} = {\left\| { - {\alpha_k^*}U({\gamma _*}){V }({\gamma _*})^\dag} \right\|_2} =\alpha_k^*{\left\| {WV({\gamma _*}){V }({\gamma _*})^\dag} \right\|_2} = \alpha_k^*,
\end{equation*}
and
\begin{equation}\label{dd}
\Delta {V }({\gamma _*})  =  - {\alpha_k^*}U({\gamma _*})\qquad \Leftrightarrow \qquad\Delta {v_i}({\gamma _*})=- {\alpha_k^*}{u_i}({\gamma _*}),\qquad i=1,\ldots,k.
\end{equation}
Using (\ref{a}) and (\ref{dd}) for every $i=1,\ldots,k,$ yields
\begin{equation}\label{bd}
\left( {{B_i} + \Delta} \right){v_i}({\gamma _*}) = - \sum\limits_{p = 1}^{k - i} {{\gamma _*}_{i,p}} {v_{i+p}}({\gamma _*}).
\end{equation}

By Lemma \ref{fullrank}, we know that $v_i(\gamma_*), (i=1,\ldots,k)$ are linearly independent vectors. Consequently, $\left( {{B_i} + \Delta} \right){v_j}({\gamma _*})$ is a nonzero vector for every $i,j \in \left\{ {1, \ldots ,k} \right\},$ except for the case $i=j=k.$ Now, define the $k$ vectors $\psi_k, \psi_{k-1},\ldots,\psi_1$ by
\begin{equation}\label{psiha}
{\psi_k} = {v_k(\gamma_*)},\qquad { \mbox{and}} \qquad {\psi_i} = \left( {\prod\limits_{p = i + 1}^k {\left( {{B_p} + \Delta } \right)} } \right){v_i(\gamma_*)},\qquad i = k - 1, \ldots, 1,
\end{equation}
note that again applying (\ref{bd}) concludes $\psi_k, \psi_{k-1},\ldots,\psi_1$ as $k$ nonzero vectors. It can be easily verified that $\left( {{B_i} + \Delta} \right)$ and $\left( {{B_j} + \Delta} \right)$ are commuting matrices, i.e.,
\begin{equation}\label{com}
\left( {{B_i} + \Delta} \right)\left( {{B_j} + \Delta} \right) = \left( {{B_j} + \Delta} \right)\left( {{B_i} + \Delta} \right),\qquad i,j = 1, \hdots ,k.
\end{equation}
Let us now, show that $A+\Delta \in \mathcal{M}_k$. First, form (\ref{a}) (but for $\gamma=\gamma_*$) and (\ref{dd}) we have
\begin{equation*}
{\left( {B_k  + \Delta } \right)\psi_k  = \left( {B_k  + \Delta } \right)v_k (\gamma _* ) =  {\alpha_k^*}{u_k}({\gamma _*})- {\alpha_k^*}{u_k}({\gamma _*})=0,}
\end{equation*}
next by considering the equations (\ref{a}), (\ref{bd}) and (\ref{com}) we can derive that
\begin{eqnarray*}
\left( {B_{k - 1}  + \Delta } \right)\psi_{k - 1}&=&\left( {B_{k - 1}  + \Delta } \right)\left( {B_k  + \Delta } \right)v_{k - 1} (\gamma _* )\\& =& \left( {B_k  + \Delta } \right)\left( {B_{k - 1}  + \Delta } \right)v_{k - 1} (\gamma _* ) \\&=& \left( {B_k  + \Delta } \right)\left( { - {\gamma _*}_{1,k - 1} v_k (\gamma _* )} \right) \\&=& 0,
\end{eqnarray*}
using the above results and performing analogous calculations, concludes
\begin{eqnarray*}
\left( {B_{k - 2}  + \Delta } \right)\psi_{k - 2}  &=& \left( {B_{k - 2}  + \Delta } \right)\left( {B_{k - 1}  + \Delta } \right)\left( {B_k  + \Delta } \right)v_{k - 2} (\gamma _* ) \\&=& \left( {B_k  + \Delta } \right)\left( {B_{k - 1}  + \Delta } \right)\left( {B_{k - 2}  + \Delta } \right)v_{k - 2} (\gamma _* )\\& =& \left( {B_k  + \Delta } \right)\left( {B_{k - 1}  + \Delta } \right)\left( { - {\gamma _*} _{1,k - 2} v_{k - 1} (\gamma _* ) - {\gamma _*} _{2,k - 2} v_k (\gamma _* )} \right) \\&=&
\left( {B_k  + \Delta } \right)\left( {B_{k - 1}  + \Delta } \right)\left( { - {\gamma _*} _{1,k - 2} v_{k - 1} (\gamma _* )} \right) \\&+& \left( {B_{k - 1}  + \Delta } \right)\left( {B_k  + \Delta } \right)\left( { - {\gamma _*} _{2,k - 2} v_k (\gamma _* )} \right)\\&=& 0.
\end{eqnarray*}
By following similar computations it is straightforwards to verify that
\begin{equation*}
{\left( {{B_i} + \Delta } \right)}\psi_i=0,\qquad\Leftrightarrow\qquad(A+\Delta)\psi_i=\lambda_i \psi_i,\qquad i = k, \ldots ,1,
\end{equation*}
This means that each $\lambda_i, (i = 1, \ldots ,k)$ is an eigenvalue of $A+\Delta$ corresponding to each $\psi_i$ as an associated eigenvector.

The main results of this section are summarized in the next theorem.
\begin{theorem}
Let $A \in \mathbb{C}^{n \times n}$ and $k\le n$ complex numbers $ \lambda _1 ,\lambda _2 , \hdots ,\lambda _k$ be given. If $\gamma_*$ is a point where the singular value ${s_\kappa }\left( {Q_A(\gamma )} \right)$ attains its maximum value, such that $\mathop \Pi \limits_{i = 1}^{k - 1} \gamma_* {_{i,1}} > 0$, and if $\alpha_k^*>0$ is a simple singular value, then $ A+\Delta \in \mathcal{M}_k$ and  ${\left\| \Delta  \right\|_2} = \alpha_k^*$, where $\Delta$ is given by (\ref{delta}).
\end{theorem}
\begin{remark}
In this section, our discussion concerned the construction of an optimal perturbation $\Delta$, benefiting the results obtained in the previous section. In particular, one of the most significant issues in Section \ref{prop}, is defining the two quantities $\gamma_*$ and $\alpha_k^*$ (see Definition \ref{sstar}) and using Theorem \ref{meng} which results in further consequences. Actually, we should follow analogous method, similar to what is described in the Sections \ref{prop} and \ref{const}, if construction of the optimal perturbation corresponding to every $\alpha_i, (i=1,\ldots, k)$ (see Corollary \ref{alfaha}) is desired. However, many numerical experiments were tried and noticed that $\alpha_k^*$ is greater than other $\alpha_i^*$ except in very rarely examples, and with an insignificant difference. Where $\alpha_i^*, (i=1,\ldots, k)$ are defined similar to Definition \ref{sstar}. Nevertheless, if $\alpha_k^*$ is considered it will do and will include the majority of cases.
\end{remark}
\section{Conclusion}\label{conc}
Lippert \cite{lipertk} in a challenging paper established a geometric motivation for spectral norm distance from $A \in \mathbb{C}^{n \times n}$ to the set $\mathcal{M}_k$ of matrices that have $k$ prescribed eigenvalues. In this work, a clear computational method for obtaining the optimal perturbation $\Delta$ such that $A+\Delta \in \mathcal{M}_k$ was presented. For the spectral norm distance from $A$ to $\mathcal{M}_k$, some lower bounds were obtained. Furthermore, an optimal perturbation to $A$, such that the perturbed matrix has the same $k$ eigenvalues, was constructed under two qualifications, that were, $s_*$ be a simple singular value and $\mathop \Pi \limits_{i = 1}^{k - 1} \gamma_* {_{i,1}} \neq 0$. We presented a numerical technique which can be considered as a generalization of results obtained in \cite{gracia,lipert,mengi}. Finally, it is pointed out that if one or both of the two assumptions are violated then the procedure described in Section \ref{const} still works, but the perturbation $\Delta$ is not necessarily optimal. In this situation, ${\left\| \Delta  \right\|_2}$ is as an upper bound for $\rho_2 (A,\mathcal{M}_k)$ whereas, we have lower bounds obtained in Section \ref{lowe} for this distance. However, finding the minimum norm perturbation $\Delta$ when $\mathop \Pi \limits_{i = 1}^{k - 1} \gamma_* {_{i,1}} = 0$, or $s_*$ has multiplicity greater than one, remains for a future investigation.
% % % % % % % % % % % % % % % % % % % % % % % % % % % % % % % % % % % % % % % % % % % % % % % % % % % % % % % % % %
\section{Applications and numerical experiment}
\label{exapmle}
In this section, we are interested to present some conceivable applications of the topic of the paper. In addition, a numerical experiment is applied to illustrate the validity of the method described in previous sections and clarify our aim. All computations were performed in MATLAB with 16 significant digits, however, for simplicity all numerical results are shown with 4 decimal places. Let us review concisely the main topic of the article. In this paper, we computed the  distance for $A\in \mathbb{C}^{n\times n}$ to $n\times n$ matrices that had $k\le n$ given complex numbers $\lambda_1,\ldots\lambda_k$ in their spectrum, while constructing an associated perturbation of $A$ was also considered. Nevertheless, another aspect of the problem may be investigated. Let us now concentrate on the subject of \textit{finding a matrix with some ordered eigenvalues}, and consider some topics and utilizations of matrices, that at the moment of this writing we find important for our aim and may be of general interest to the readers as well.

\textbf{Inverse eigenvalue problem.}
One of straightforward usage of aforesaid viewpoint is computing a matrix that has some prespecified eigenvalues, which is known as \textit{inverse eigenvalue problem}. An inverse eigenvalue problem (IEP) concerns the reconstruction of a matrix from prescribed spectral data. The spectral data involved may consist of the complete or only partial information of eigenvalues or eigenvectors. Inverse eigenvalue problem has a long list of applications in areas such as control theory, mechanics, signal processing and numerical analysis. Two acceptable references for theory and application of inverse eigenvalue problem are \cite{chu, chugolub}. Assume now, we are asked to find a matrix having given scalars $\lambda_1,\ldots, \lambda_l \in \mathbb{C}$ where $l \le n$, as its eigenvalues. To do this, by using the technique explained in this paper, one can consider an arbitrary matrix, namely, $A$ in the craved size. Next, by following procedure described in Section \ref{const}, the desired matrix which has $\lambda_1,\ldots,\lambda_l$ as some of its eigenvalues is computable.

\textbf{Approximating a matrix with another one that has prescribed eigenvalues.} In many applications of matrices, particularly in optimization problems, we are concerned with a matrix that must have some assumptions on its eigenvalues. Suppose, for instance, that $f(x)$ is a function that we want to maximize or minimize it. In the various method for optimizing $f(x)$, it is necessary that $\nabla^2f$ be a positive definite matrix. A matrix $M$ is positive definite if and only if its all eigenvalues are greater than zero. However, $\nabla^2f$ may not be positive definite for some cases. In such problems, some approaches will be applied to overcome this difficulty which all of them benefit the constructing a correction matrix $\Delta M$ such that ensures all eigenvalues of $M+\Delta M$ are greater enough than zero. For a symmetric matrix $M$, we can obtain the correction matrix $\Delta M$ with minimum Euclidean norm and Frobenius norm by using the techniques performed in Section 3.4 of \cite{nocedal}. Now, our procedure for a general non-positive definite matrix $M$, can be applied. First, construct the matrix $\Delta M$ and then replace $M$ by a positive definite approximation $M+\Delta M$, in which all negative eigenvalues in $M$ are replaced by arbitrary (large enough) positive numbers, while keeping its other (positive) eigenvalues unchanged.
\begin{example}
Consider the non-positive definite matrix
\begin{equation*}
A = \left[ {\begin{array}{*{20}c}
   3 & 6 & 9 & {10}  \\
   4 & 1 & { - 1} & { - 2}  \\
   7 & 5 & 0 & { - 4}  \\
   4 & { - 3} & { - 1} & 6  \\
\end{array}} \right],
\end{equation*}
which is a random matrix generated by MATLAB and its eigenvalues are 12.9377, 7.0550, -0.3641 and -9.6286. Assume now the set $\Lambda=\left\lbrace 12.9377, 7.0550, \varepsilon, \varepsilon\right\rbrace $ in order to construct a positive definite approximation of $A$ where $\varepsilon=1\times 10^{-4}$ is chosen arbitrarily. Using the MATLAB function \texttt{fmincon} which finds a constrained minimum of a function of several variables, one can derive that the function $s_{13}(Q_A(\gamma))$ attains its maximum value at the finite point $\gamma_{*}=[6.6686, 8.2009, 7.9580, 0.5925-0.0021 i, -2.3778-0.0009 i, 1.4475-0.0009 i]$ and  $\alpha_4^*=s_{13}(Q_A(\gamma_*))=5.1231,$ is an isolated singular value of constant matrix $Q_A(\gamma_*)$. Furthermore, $\left\| {U(\gamma _* )^* U(\gamma _* ) - V(\gamma _* )^* V(\gamma _* )} \right\| = {{2}}{{.0266}}\times 10^{-4}$ verifying Corollary \ref{natijeh}. By applying the technique explained in the paper, the matrix $\Delta$ in (\ref{delta}) is computed as follows:
\[
\Delta  = \left[ {\begin{array}{*{20}c}
   {{ {3}}{ {.4650  -  0}}{ {.0011i}}} & {{ { - 1}}{ {.1376  -  0}}{ {.0004i}}} & {{ {1}}{ {.2052  -  0}}{ {.0006i}}} & {{ {3}}{ {.3916  -  0}}{ {.0025i}}}  \\
   {{ {0}}{ {.2301  -  0}}{ {.0022i}}} & {{ {1}}{ {.4468  -  0}}{ {.0012i}}} & {{ { - 4}}{ {.5382  -  0}}{ {.0019i}}} & {{ {1}}{ {.8719  -  0}}{ {.0046i}}}  \\
   {{ { - 0}}{ {.2479  +  0}}{ {.0006i}}} & {{ {4}}{ {.6130  -  0}}{ {.0005i}}} & {{ {1}}{ {.9171  +  0}}{ {.0005i}}} & {{ {1}}{ {.1181  +  0}}{ {.0035i}}}  \\
   {{ { - 3}}{ {.7579  +  0}}{ {.0001i}}} & {{ { - 1}}{ {.2664  -  0}}{ {.0010i}}} & {{ {0}}{ {.7100  +  0}}{ {.0005i}}} & {{ {3}}{ {.1641  +  0}}{ {.0018i}}}  \\
\end{array}} \right],
\]
which satisfies $\left\| \Delta\right\|\simeq \alpha_4^*$ and $12.9377, 7.0550, 10^{-4}, 10^{-4}$ are eigenvalues of $A+\Delta$. Clearly, $A+\Delta$ is a positive definite matrix. 
\end{example}
Moreover, this idea can be implemented in areas such as solving a singular or ill condition system of linear equations, i.e., $Ax=b$ for some $A\in \mathbb{C}^{n\times n}$, and studying the solutions of fast-slow systems. These two issues are reviewed briefly in the following.

Notice that, singular square matrices are a thin subset of the space of all square matrices and adding a tiny random perturbation to a singular matrix makes it nonsingular. But the perturbed matrix, namely, $\tilde{A}$ may not be useful in practice. In other words, we will able to compute a solution for $\tilde{A}x=b$, but that solution may be meaningless. This failure, in general, is caused mainly by changing all eigenvalues of $A$ and being large enough the condition number of $\tilde{A}$ in which the smallest eigenvalue of $\tilde{A}$ has a significant effect. It seems that we can cope simultaneously with these failures by applying the method explained in Section \ref{const}. For the case of a singular (ill condition) system we construct a perturbation in which all zero (small enough) eigenvalues in $A$ are replaced by nonzero (large enough) numbers. Nevertheless, it should be noted that the accuracy and credibility of the solution depends extremely on the problem we are trying to solve. This topic is known as \textit{regularization}, however, it may have different names in different contexts.

On the other hand, consider the fast-slow systems. In order to apply some methods for studying the solutions of fast-slow systems, we must have a matrix that has no eigenvalues on the imaginary axis \cite{feni}.We can overcome this problem, by constructing a perturbation of the matrix such that all its eigenvalues are the same and changing only the eigenvalues which are on the imaginary axis. See \cite{hopp,gras,chiba} and references therein for theory and applications of fast-slow systems.

\textbf{Low-rank approximation problem.} Suppose that $A \in \mathbb{C}^{m\times n}$ and a positive integer $r<\text{min}\{m,n\}$ are given. In mathematics, \textit{low-rank approximation problem} concerns finding a matrix $D\in \mathbb{C}^{m\times n}$ such that rank$(D)=r$ and $D$ is nearest matrix to $A$, i.e., $\lVert A-D\rVert\ge 0$ is as small as possible and $\|.\|$ can in general be any matrix norm.  The problem is used for mathematical modeling, data and image compression, noise reduction, system and control theory, computer algebra and so on.  Many different techniques are provided for computing a matrix with another one of the same size that has reduced rank. See \cite{chulow, mark1,mark2,mark3} as the suggested references on the low-rank approximation problem and its applications. However, we can use the results obtained in this paper for computing a low-rank approximation of a given square matrix. To do this, by applying the procedure described, we can set $\lambda_1=\ldots=\lambda_m=0$ where $m\le k.$ Now if its corresponding eigenvectors $\psi_1,\ldots,\psi_m$ introduced in (\ref{psiha}) satisfy rank$\left( \left[\psi_1,\ldots,\psi_m \right] \right) =\tau$, then it is clear that rank of the perturbed matrix $A+\Delta$ is almost $n-\tau$.

%%%%%%%%%%%%%%%%%%%%%%%%%%%%%%%%%%%%%%%%%%%%%%%%%%%%%%%%%%%%%%%%%%%%%%%%%%%%%%%%%%%%%%%%%%%%%%%%%%%%%%%%%%%%%%%%%%%%%%%%%%%%%%%%%%%%%%%%%%%%%%%%
\section{Acknowledgment}
The authors are very grateful to Prof. Ross A. Lippert who provided constructive comments and valuable suggestions which have considerably improved the article.

%\pagebreak
%%%%%%%%%%%%%%%%%%%%%%%%%%%%%%%%%%%%%%%%%%%%%%%%%%%%%%%%%%%%%%%%%%%

%%%%%%%%%%%%%%%%%%%%%%%%%%%%%%%%%%%%%%%%%%%%%%%%%%%%%%%%
\end{document}